\documentclass[10pt]{article}

\usepackage{amsthm}
\usepackage{amsmath}
\usepackage{amsfonts}
\usepackage{amssymb}
\usepackage{mathrsfs}
\usepackage{empheq}
\usepackage{stmaryrd}
\usepackage{dsfont}
\usepackage{enumerate}
\usepackage{comment}
\usepackage{pgfplots,soul}
\usepackage[applemac]{inputenc}
\usepackage[colorlinks=true,urlcolor=blue]{hyperref}

\newcommand{\R}{\mathbb{R}}
\newcommand{\C}{\mathbb{C}}
\newcommand{\K}{\mathbb{K}}
\newcommand{\vlL}{\mathbb{L}}

\newcommand{\Dr}{\mathscr{D}}

\newcommand{\Lr}{\mathscr{L}}

\newcommand{\vphi}{\varphi}

\newcommand{\dsp}{\displaystyle}
\newcommand{\ovl}{\overline}
\newcommand{\udl}{\underline}

\newcommand{\vint}{\int\limits}

\newcommand{\inj}{\hookrightarrow}
\newcommand{\tends}{\longrightarrow}

\newcommand{\loc}{\mathrm{loc}}

\renewcommand{\d}{\mathrm{d}}
\newcommand{\vo}{\mathrm{o}}

\newcommand{\vi}{\mathrm{i}}

\renewcommand{\le}{\leqslant}
\renewcommand{\ge}{\geqslant}
\renewcommand{\Re}{\mathrm{Re}}
\renewcommand{\Im}{\mathrm{Im}}

\newcommand{\p}{\prime}

\newcommand{\eqdef}{\stackrel{\mathrm{def}}{=}}

\numberwithin{equation}{section}

\newtheorem{thm}{Theorem}[section]

\theoremstyle{definition}

\newtheorem{defi}[thm]{Definition}

\newenvironment{proof*}{\noindent{\bf Proof.}}{\qed}

\title{\huge \sc The dual space of a complex Banach space restricted to the field of real numbers}
\author{\sc Pascal Bégout}
\date{}

\begin{document}

\maketitle

\begin{center}
Institut de Mathématiques de Toulouse \& TSE	\\
Université Toulouse I Capitole				\\
1, Esplanade de l’Université				\\
31080 Toulouse Cedex 6, FRANCE
\bigskip \\
\text{
{\footnotesize E-mail\:: \href{mailto:Pascal.Begout@math.cnrs.fr}{\udl{\texttt{Pascal.Begout@math.cnrs.fr}}}}
}
\end{center}

\begin{abstract}
Solutions of some partial differential equations are obtained as critical points of a real funtional. Then the Banach space where this functional is defined has to be real, otherwise, it is not differentiable. It follows that the equation is solved with respect to the real dual space of this Banach space. But if the solution is complex-valued there is the following problem: what does the multiplication of this equation by a complex number mean~? In this note, we explain how to rigorously define this operation.
\end{abstract}

{\let\thefootnote\relax}
{\let\thefootnote\relax}

\tableofcontents

\baselineskip .6cm

\section{Introduction}
\label{secint}

Some partial differential equations are obained in the following way. Let $X$ be a real Banach space and let $F\in C^1(X;\R).$ It follows that $F^\p\in C(X;X^\star).$ If $u$ is a critical point of $F$ then $u$ solves,
\begin{gather}
\label{pdem}
F^\p(u)=0, \text{ in } X^\star.
\end{gather}
If $X$ is a complex Banach space then $F$ is no more differentiable and the above equation makes no sense. Hence the necessity to deal with real Banach spaces. On the other hand, if the above equation is complex-valued, we may be interested in multiplying it by $\vi.$ But $X^\star=\Lr(X;\R)$ so that $\vi F(u)\notin X^\star$ and this operation is not allowed.
\medskip
\\
Let us illustrate the problem with a very simple example. The linear space $L^2(\R^N;\C)$ may be equipped with the field of complex numbers $\C.$ The natural inner product is,
\begin{gather*}
\forall(u,v)\in L^2(\R^N;\C)\times L^2(\R^N;\C), \; ((u,v))_{L^2(\R^N;\C)}=\vint_{\R^N}u(x)\ovl{v(x)}\d x,
\end{gather*}
with the norm,
\begin{gather*}
\forall u\in L^2(\R^N;\C), \; \|u\|_{L^2(\R^N,\C)}^2=((u,u))_{L^2(\R^N;\C)}=\vint_{\R^N}|u(x)|^2\d x.
\end{gather*}
Let us denote this topological vector space by $L^2(\R^N;\C)_\C.$ Another choice is to equipp $L^2(\R^N;\C)$ with the field of real numbers $\R.$ The natural scalar product is then,
\begin{gather*}
\forall(u,v)\in L^2(\R^N;\C)\times L^2(\R^N;\C), \; (u,v)_{L^2(\R^N;\C)}=\Re\vint_{\R^N}u(x)\ovl{v(x)}\d x,
\end{gather*}
with the same norm as above. Let us denote this topological vector space by $L^2(\R^N;\C)_\R.$ In both cases, $L^2(\R^N;\C)_\C$ and $L^2(\R^N;\C)_\R$ are, respectively, complex and real Hilbert spaces with the same topology and have the same elements. But there is a fundamental difference concerning the differentiable functions. We first recall that if a function $F:A\tends B,$ between two Banach spaces $A$ and $B,$ is Fréchet-differentiable at $u\in A$ then $F^\p(u)\in\Lr(A;B),$ the space of linear continuous functions from $A$ to $B.$ In particular, $A,$ $B$ and $\Lr(A;B)$ deal with the same field. Let for any $u\in L^2(\R^N;\C),$ $F(u)=\frac12\|u\|_{L^2(\R^N;\C)}^2.$ Then for any $u\in L^2(\R^N;\C)$ and $v\in L^2(\R^N;\C),$
\begin{gather*}
F(u+v)=F(u)+\Re\vint_{\R^N}u(x)\ovl{v(x)}\d x+\frac12\|v\|_{L^2(\R^N;\C)}^2,
\end{gather*}
and $\frac12\|v\|_{L^2(\R^N;\C)}^2=\vo\left(\|v\|_{L^2(\R^N;\C)}\right).$ It follows that for any $u\in L^2(\R^N;\C),$ the function
\begin{gather}
\label{linfor}
f_u:v\longmapsto\Re\vint_{\R^N}u(x)\ovl{v(x)}\d x,
\end{gather}
is $\R$-linear and continuous from $L^2(\R^N;\C)_\R$ to $\R.$ So $F$ is differentiable over $L^2(\R^N;\C)_\R$ and for any $u\in L^2(\R^N;\C),$ $F^\p(u)=f_u.$ But $f_u$ is not $\C$-linear and furthermore $F:L^2(\R^N;\C)_\C\tends\R$ is never differentiable (except for $u=0).$ More generally, if $X_\C$ is a complex Banach space and if $F:X_\C\tends\R$ then $F$ is never differentiable at $u\in X_\C$ (except for the constant functions) because if it was, we should have for any $v\in X_\C,$
\begin{gather*}
\langle F^\p(u),\vi v\rangle\in\R,	\\
\langle F^\p(u),\vi v\rangle=\vi\langle F^\p(u),v\rangle\in\vi\R.
\end{gather*}
Hence, $F^\p(u)=0.$ It follows that, as soon as a differential of a real function is needed, the functional space where it lies has to be considered over the field $\R.$
\medskip
\\
Now, let us consider a partial differential equation whose solutions are complex-valued. For most of these equations, it is required that some real functions have a differential (or a derivative), such as, for instance,  an energy,
\begin{gather*}
E(u)=\frac12\int_\Omega|\nabla u(x)|^2\d x+\frac1p\int_\Omega|u(x)|^p\d x,
\end{gather*}
for some $p\ge1.$ Some of these solutions are also obtained by minimizing a real functional. See Section~\ref{secexa} for some examples. From the discussion above, it is clear that the Banach spaces considered have to be over $\R.$ But then, there is another problem. Let us illustrate it by the example of the linear Schrödinger equation,
\begin{gather*}
\vi\frac{\partial u}{\partial t}+\Delta u=0, \quad (t,x)\in\R\times\R^N.
\end{gather*}
If we endow $L^2(\R^N;\C)$ of the field of real numbers and take an initial value in this space, it is well-known that the solution $u$ belongs to
\begin{gather*}
C\big(\R;L^2(\R^N;\C)\big)\cap C^1\big(\R;H^{-2}(\R^N;\C)\big),
\end{gather*}
where $H^{-2}(\R^N;\C)$ is the topological dual space of $H^2(\R^N;\C),$ which means that,
\begin{gather*}
\frac{\partial u}{\partial t}(t)\in H^{-2}(\R^N;\C)\eqdef\Lr\big(H^2(\R^N;\C);\R\big),
\end{gather*}
for any $t\in\R.$ But then, what does the product $\vi\dfrac{\partial u}{\partial t}$ mean\,? This question can be answered by looking at the method used to solve the equation, or by the equation itself. But from a mathematical point of view, this method is not satisfying because it has been adapted to each equation and the answer is not exogenous. More generally, if $X_\R$ is a real Banach space which is the restriction to the field $\R$ of a complex Banach space $X_\C,$ and if $T\in X_\R^\star\eqdef\Lr(X_\R;\R),$ may we give a meaning to the product $\vi T\in\Lr(X_\R;\R)\,?$ An answer may be given if, additionally, $X_\C$ is a complex Hilbert space whose inner product is $((\:.\:,\:.\:)).$ Then $X_\R$ is a real Hilbert space whose scalar product is $(\:.\:,\:.\:)\eqdef\Re((\:.\:,\:.\:)).$ We have the natural injection $X_\R\overset{T}{\inj}X_\R^\star$ where for $u\in X_\R,$ $T_u\in X_\R^\star$ is defined by,
\begin{gather}
\label{RFT}
\forall v\in X_\R, \; \langle T_u,v\rangle_{X_\R^\star,X_\R}=(u,v).
\end{gather}
But if $v\in X_\R,$ so is $\vi v$ and for any $v\in X_\R,$ $\langle T_u,\vi v\rangle_{X_\R^\star,X_\R}=-(\vi u,v).$ So it is natural to define $\vi T\in X^\star_\R$ by,
\begin{gather}
\label{defiTH}
\forall v\in X_\R, \; \langle\vi T_u,v\rangle_{X_\R^\star,X_\R}=(\vi u,v)=-\langle T_u,\vi v\rangle_{X_\R^\star,X_\R}.
\end{gather}
In particular, when an equation makes sense in $X^\star_\R,$ \eqref{defiTH} is a way to give a meaning to the sentence ``We multiply the equation by $\vi$''. By the Riesz-Fréchet representation Theorem, it follows that the map $T:X_\R\tends X_\R^\star$ defined by~\eqref{RFT} is an isometric isomorphism and it completely answers to the question. But this method fails if $X_\C$ is not Hilbert space. Formula~\eqref{defiTH} suggests a relation between the real dual space of $X_\R$ and the (complex) adjoint space of $X_\C,$ as we shall see in Section~\ref{secdef}.
\medskip
\\
Unfortunately, some authors deal with complex Banach space $X$ and solve some equations as~\eqref{pdem}, with $F:X\tends\R,$ which has no sense since $F$ is never differentiable. The goal of this note is to provide a rigorous method to avoid this kind of mistake.
\medskip
\\
This paper is organized as follows. In the following section, we give some notations which will be used throughout this paper, and we recall some well-known facts on functional analysis. In Section~\ref{secdef}, we give a simple result which allow to give a natural meaning to the product $\vi T,$ for $T\in\Lr(X;\R),$ in order to have $\vi T\in\Lr(X;\R)$ (Definition~\ref{defiT}). In Section~\ref{secapp}, we give some applications of Definition~\ref{defiT}. In the last section, we give some example where Definition~\ref{defiT} is implicitely used by some authors.

\section{Notations}
\label{secnot}

Let $X$ be a $\C$-linear space. If $X$ is a $\C$-topological vector space, it will be denoted by $X_\C$ and its restriction to the field $\R$ will be denoted by $X_\R.$ It follows that $X_\R$ is a $\R$-topological vector space and has the same topology and the same elements than $X_\C.$ This algebraic set will be simply denoted by $X.$ This notation will be used when no consideration of vector space structure and topolopy will be needed. If $X_\C$ is a normed space then $X_\R$ has the same norm which will be denoted by $\|\:.\:\|_X$ for the both spaces. If $H_\C$ is a complex Hilbert space, its inner product will be denoted by $((\:.\:,\:.\:))_{H_\C}.$ It follows that $H_\R$ is a real Hilbert space whose the scalar product is denoted by $(\:.\:,\:.\:)_{H_\R}\eqdef\Re((\:.\:,\:.\:))_{H_\C}.$ Let $\K\in\{\R,\C\}.$ For $X_\K$ a $\K$-topological vector space, its $\K$-topological dual space $\Lr(X_\K;\K)$ will be denoted by $X^\star_\K,$ and its $X^\star_\K-X_\K$ duality product will be written $\langle\: . \; , \: . \:\rangle_{X^\star_\K,X_\K}\in\K.$ It is well-known that the $\K$-topological dual space $X_\K^\star$ is complete and if $X_\K$ is normed then $X_\K^\star$ is a $\K$-Banach space. For a $\C$-topological vector space $X_\C,$ a function $f:X_\C\tends\C$ is said to be an \textit{anti-linear form} (or a \textit{semilinear form}) if $f$ is additive and if for any $(\lambda,x)\in\C\times X_\C,$ $f(\lambda x)=\ovl\lambda f(x),$ where $\ovl\lambda$ is the conjugate complex number of $\lambda.$ The \textit{anti-dual space} or \textit{adjoint space} of $X_\C$ is denoted by $\ovl{X^\star}.$ It consists of the continuous anti-linear forms on $X_\C.$ It is a $\C$-complete vector space. Its $\ovl{X^\star}-X_\C$ duality product will be written $\langle\: . \; , \: . \:\rangle_{\ovl{X^\star},X_\C}\in\C.$ If $X_\C$ is a normed space then $\ovl{X^\star}$ is a $\C$-Banach space. Note that $X^\star_\C,$ $\ovl{X^\star}$ and $X^\star_\R$ have the same norm (Kato~\cite{MR1335452}: p.10-14, p.134) and will be simply denoted by $\|\:.\:\|_{X^\star}$ for these three spaces. Let $X_\K$ and $Y_\K$ be locally convex Hausdorff $\K$-topological vector spaces. The notation $X_\K\overset{e}{\inj}Y_\K$ means that $e:X_\K\tends Y_\K$ is $\K$-linear, continuous and one-to-one. We recall that if this embedding is dense then $Y_\K^\star\overset{e^\star}{\inj}X_\K^\star,$ where $e^\star$ is the transpose of $e:$
\begin{gather*}
\forall L\in Y_\K^\star, \; \forall x\in X, \; \langle e^\star(L),x\rangle_{X_\K^\star,X_\K}\eqdef\langle L,e(x)\rangle_{Y_\K^\star,Y_\K}.
\end{gather*}
If, furthermore, $X_\K$ is reflexive then the embedding $Y_\K^\star\overset{e^\star}{\inj}X_\K^\star$ is dense. Often, $e$ is the identity function, so that $e^\star$ is nothing else but the restriction to $X_\K$ of continuous $\K$-linear forms on $Y_\K.$ For more details, see Trèves~\cite[Corollary~5, p.188; Corollary, p.199; Theorem~18.1, p.184]{MR2296978}. Finally, for a complex number $z,$ $\Re(z)$ and $\Im(z)$ are the real part and the imaginary part of $z,$ respectively.

\section{Definition of the product of an element of a real dual space by a complex number}
\label{secdef}

In this section, we give a meaning to the product $\vi T\Lr(X;\R),$ when $T\in\Lr(X;\R).$ The following result is a trivial adaptation of the result of Brezis~\cite[Proposition~11.22, p.361]{MR2759829}, and so we omit the proof.

\begin{thm}
\label{thmmain}
Let $X_\C$ be a complex topological vector space. Consider the following map,
\begin{gather}
\label{L}
\begin{array}{rcl}
\vlL:\ovl{X^\star}		&	     \tends		&	X^\star_\R,	\medskip \\
	    	       T		&	\longmapsto	&	\Re\:T.
\end{array}
\end{gather}
Then $\vlL$ is $\R$-linear, bijective and for any $T\in X^\star_\R,$
\begin{gather*}
\forall x\in X, \; \langle\vlL^{-1}(T),x\rangle_{\ovl{X^\star},X_\C}=\langle T,x\rangle_{X^\star_\R,X_\R}+\vi\langle T,\vi x\rangle_{X^\star_\R,X_\R}.
\end{gather*}
If, in addition, $X_\C$ is a normed space then $\vlL$ is an isometry.
\end{thm}

\begin{defi}
\label{defiT}
Let $X_\C$ be a complex topological vector space and let $\vlL$ be given by Theorem~\ref{thmmain}. For any $T\in X^\star_\R,$ we define $\vi T\in X^\star_\R$ by,
\begin{gather}
\label{defiT1}
\vi T=\vlL\left(\vi\vlL^{-1}(T)\right).
\end{gather}
It follows that,
\begin{gather}
\label{iTgen}
\langle\vi T,x\rangle_{X^\star_\R,X_\R}=\langle T,-\vi x\rangle_{X^\star_\R,X_\R},
\end{gather}
for any $x\in X.$ In other words, the sentence ``We multiply $T$ by $\vi$'' has to be understood as ``We take the $X^\star_\R-X_\R$ duality product of $T$ with $-\vi x,$ for $x\in X$''.
\end{defi}

\section{Applications}
\label{secapp}

In this section, we give some applications of the previous section to give a meaning to $\vi T\in X^\star_\R,$ when $T\in X^\star_\R.$

\subsection*{The Hilbert spaces}

Let $H_\C$ be a complex Hilbert space. By Theorem~\ref{thmmain}, \eqref{iTgen} and the Riesz-Fréchet representation theorem it follows that for any $L\in H^\star_\R,$ there exists a unique $u_L\in H$ such that
\begin{gather}
\label{h1}
\langle L,v\rangle_{H^\star_\R,H_\R}=(u_L,v)_{H_\R},		\medskip \\
\label{h2}
\langle\vi L,v\rangle_{H^\star_\R,H_\R}=\langle L,-\vi v\rangle_{H^\star_\R,H_\R}=(u_L,-\vi v)_{H_\R},
\end{gather}
for any $v\in H.$ In addition, $\|u_L\|_X=\|L\|_{X^\star}=\|\vi L\|_{X^\star}.$ Let $u\in H.$ Set for any $v\in H,$ $\langle L_u,v\rangle_{H^\star_\R,H_\R}=(u,v)_{H_\R}.$ Then $L_u\in H^\star_\R$ and \eqref{h1}--\eqref{h2} are nothing else but,
\begin{gather}
\label{h3}
(\vi u,v)_{H_\R}=(u,-\vi v)_{H_\R},
\end{gather}
for any $v\in H.$ From this point of view, the formula~\eqref{iTgen} is an extension of~\eqref{h3} to the topological vector spaces and~\eqref{iTgen} reads as $\vi L_u=L_{\vi u}.$

\subsection*{The Lebesgue spaces}

Let $1\le p<\infty.$ Let $1<p^\p\le\infty$ its conjugate defined by, $\frac1p+\frac1{p^\p}=1,$ and let $\Omega$ be a nonempty open subset of $\R^N.$ It is well-known that for any $T\in\ovl{L^p(\Omega;\C)^\star},$ there exists a unique $u\in L^{p^\p}(\Omega;\C)$ such that,
\begin{gather*}
\langle T,v\rangle_{\ovl{L^p(\Omega;\C)^\star},L^p(\Omega;\C)_\C}=\vint_\Omega u(x)\ovl{v(x)}\d x,
\end{gather*}
for any $v\in L^p(\Omega;\C).$ In addition, $\|T\|_{\ovl{L^p(\Omega;\C)^\star}}=\|u\|_{L^{p^\p}(\Omega;\C)}$ (Kato~\cite[Example~1.25, p.135]{MR1335452}). Applying Theorem~\ref{thmmain} and \eqref{iTgen}, we get that the following map,
\begin{gather*}
\begin{array}{rcll}
\Phi_p:L^{p^\p}(\Omega;\C)_\R	&	     \tends		&	\big(L^p(\Omega;\C)\big)^\star_\R	&	\medskip \\
					     u	&	\longmapsto	&	\Phi_p(u):L^p(\Omega;\C)_\R		&	\tends\R	\medskip \\
						&				&	\quad\qquad\qquad\qquad	 v		&	\longmapsto\Re\dsp\vint_\Omega u(x)\ovl{v(x)}\d x.
\end{array}
\end{gather*}
is an isometric isomorphism. In other words, if $T\in L^p(\Omega;\C)^\star_\R,$ there exists a unique $u\in L^{p^\p}(\Omega;\C)$ such that,
\begin{gather*}
\langle T,v\rangle_{L^p(\Omega;\C)^\star_\R,L^p(\Omega;\C)_\R}=\Re\vint_\Omega u(x)\ovl{v(x)}\d x,	\\
\langle\vi T,v\rangle_{L^p(\Omega;\C)^\star_\R,L^p(\Omega;\C)_\R}=\Re\left(\vi\vint_\Omega u(x)\ovl{v(x)}\d x\right),
\end{gather*}
for any $v\in L^p(\Omega;\C),$ where $u=\Phi_p^{-1}(T),$ and $\vi T=\Phi_p(\vi u).$ Furthermore, $\|u\|_{L^{p^\p}(\Omega;\C)}=\|T\|_{L^p(\Omega;\C)^\star}=\|\vi T\|_{L^p(\Omega;\C)^\star}.$ We recover that,
\begin{gather*}
L^p(\Omega;\C)^\star_\R\overset{\Phi_p^{-1}}{\cong}L^{p^\p}(\Omega;\C)_\R,
\end{gather*}
and for any $u\in L^{p^\p}(\Omega;\C),$
\begin{gather*}
\langle u,v\rangle_{L^{p^\p}(\Omega;\C)_\R,L^p(\Omega;\C)_\R}=\Re\vint_\Omega u(x)\ovl{v(x)}\d x,
\end{gather*}
for any $v\in L^p(\Omega;\C).$ Note that for $p=2,$ this means that the Riesz-Fréchet and the Riesz Theorems coincide.

\subsection*{The space of distributions}

We first recall that the topological dual space $\Dr(\Omega;\C)^\star_\R$ of $\Dr(\Omega;\C)_\R$ is the space of distributions $\Dr^\p(\Omega;\R):$ weak$\star$ and strong topologies are the same on bounded sets (Schwartz~\cite[Théorèmes~XIII and XIV, p.74-75]{MR0209834}). Now, let us define $\Phi:L^1_\loc(\Omega;\C)_\C\tends\ovl{\Dr(\Omega;\C)^\star}$ by,
\begin{gather*}
\forall f\in L^1_\loc(\Omega;\C), \; \langle\Phi_f,\vphi\rangle_{\ovl{\Dr(\Omega;\C)^\star},\Dr(\Omega;\C)_\C}=\vint_\Omega f(x)\ovl{\vphi(x)}\d x,
\end{gather*}
for any $\vphi\in\Dr(\Omega;\C).$ Standard arguments show that $\Phi$ is $\C$-linear, continuous and one-to-one. Applying Theorem~\ref{thmmain} and \eqref{iTgen}, if follows that,
\begin{gather*}
L^1_\loc(\Omega;\C)_\R\overset{T}{\inj}\Dr^\p(\Omega;\R),
\end{gather*}
where $T=\Re(\Phi).$ We then have for any $f\in L^1_\loc(\Omega;\C)_\R,$ $\vi f\in L^1_\loc(\Omega;\C)_\R$ and,
\begin{gather*}
\langle T_f,\vphi\rangle_{\Dr^\p(\Omega;\R),\Dr(\Omega;\C)_\R}=\Re\vint_\Omega f(x)\ovl{\vphi(x)}\d x,	\\
\langle\vi T_f,\vphi\rangle_{\Dr^\p(\Omega;\R),\Dr(\Omega;\C)_\R}
=\langle T_{\vi f},\vphi\rangle_{\Dr^\p(\Omega;\R),\Dr(\Omega;\C)_\R}=\langle T_f,-\vi\vphi\rangle_{\Dr^\p(\Omega;\R),\Dr(\Omega;\C)_\R},
\end{gather*}
for any $\vphi\in\Dr(\Omega;\C).$

\section{Examples}
\label{secexa}

In this section, we give some examples of situations where the method used requires to restrict some complex Banach spaces to the field of real numbers.

\subsection*{The nonlinear Schrödinger equation}

To solve the following nonlinear Schrödinger equation,
\begin{gather*}
\begin{cases}
\vi\dfrac{\partial u}{\partial t}+\Delta u+g(u)=0,	&			\\
u_{|\partial\Omega}=0,					&	\dfrac{}{}	\\
u(0)=u_0,
\end{cases}
\end{gather*}
with $u_0\in H^1_0(\Omega;\C),$ some assumptions are done about $g.$ Among them, there is the existence of a $G\in C^1(H^1_0(\Omega;\C);\R)$ such that $G^\p=g.$ Such a function cannot exist if $H^1_0(\Omega;\C)$ is over the field $\C.$ The starting point is an abstract result which deals with a complex Hilbert space but with the real scalar product (Cazenave~\cite[Theorem~3.3.1, p.63-64]{MR2002047}). This theory applies, at least, where $g$ is as follows:
\begin{gather*}
g(u)=Vu+\lambda|u|^\alpha u+(W\star|u|^2)u,
\end{gather*}
with,

$\bullet$
$V\in L^p(\R^N;\R)+ L^\infty(\R^N;\R),$ $p\ge1,$ $p>\frac{N}2,$

$\bullet$
$\lambda\in\R,$ $\alpha\ge0,$ $(N-2)\alpha<4,$

$\bullet$
$W\in L^q(\R^N;\R)+ L^\infty(\R^N;\R),$ $q\ge1,$ $q>\frac{N}4,$ $W$ is even.
\\
It follows that $\frac{\partial u}{\partial t}\in H^{-1}(\R^N)\eqdef\Lr\left((H^1(\R^N;\C);\R\right)$ and the product $\vi\frac{\partial u}{\partial t}$ is understood as in Definition~\ref{defiT}. For more details, see Cazenave~\cite[Chapters~3-4]{MR2002047}.

\subsection*{The damped nonlinear Schrödinger equation}

In \cite{MR4098330,MR4053613,MR4340780,MR4503241}, the authors study the finite time extinction property of the solutions to,
\begin{gather*}
\begin{cases}
\vi\dfrac{\partial u}{\partial t}+\Delta u+V(x)u+a|u|^{-(1-m)}u=f(t,x),	&	\text{ in } (0,\infty)\times\Omega,				\\
u_{|\partial\Omega}=0,									&	\text{ on } (0,\infty)\times\partial\Omega,	\dfrac{}{}	\\
u(0)= u_0,												&	\text{ in } \Omega,
\end{cases}
\end{gather*}
where $0\le m\le1,$ $a\in\C,$ $\Omega\subseteq\R^N,$ $f\in L^1_\loc\big([0,\infty);L^2(\Omega;\C)\big)$ and $V\in L^1_\loc(\Omega;\R).$ The well posedness of the Cauchy problem is established for $u_0\in L^2(\Omega;\C),$ $u_0\in H^1_0(\Omega;\C)$ and $u_0\in H^2(\Omega;\C).$ The proof of the existence strongly relies on the theory of maximal monotone operators where the Hilbert spaces (and, more generally, the Banach spaces) are assumed to be real (Barbu~\cite{MR0390843,MR2582280}, Brezis~\cite{MR0348562}, Vrabie~\cite{MR1375237}). In particular, if $u_0\in L^2(\Omega)$ then $\frac{\partial u}{\partial t}\in H^{-2}(\Omega)+L^\frac2m(\Omega;\C),$ where $H^{-2}(\Omega)\eqdef\Lr\left(H^2_0(\Omega;\C);\R\right),$ and the product $\vi\frac{\partial u}{\partial t}$ is understood as in Definition~\ref{defiT}.

\subsection*{The stationary Schrödinger equation with periodic magnetic potential}

In \cite{MR4228668}, the authors show the existence of a solution to the following stationary Schrödinger equation with periodic magnetic potential,
\begin{gather*}
\begin{cases}
-\Delta_Au+V(x)u=\lambda g(u) \; \text{ in } \R^N,	\medskip \\
u\in H^1_{A,V}(\R^N),
\end{cases}
\end{gather*}
where $-\Delta_Au=(-\vi\,\nabla + A)^2u,$ $N\ge2,$ $\lambda>0,$ and $A\in L^N_\loc(\R^N;\R^N)$ and $V\in L^\frac N2_\loc(\R^N;\R)$ satisfy some suitable assumptions. The Hilbert space $H^1_{A,V}(\R^N;\C)$ is,
\begin{gather*}
H^1_{A,V}(\R^N;\C)=\Big\{u\in H^1(\R^N;\C); Vu^2\in L^1(\R^N;\C) \text{ and } Au\in L^2(\R^N;\C^N)\Big\},
\end{gather*}
whose the scalar product is,
\begin{gather*}
\forall u,v\in H^1_{A,V}(\R^N), \;
(u,v)_{H^1_{A,V}(\R^N)}=\Re\vint_{\R^N}Vu\,\ovl v\,\d x+\Re\vint_{\R^N}\nabla_Au.\ovl{\nabla_Av}\,\d x,
\end{gather*}
where $\nabla_A=\nabla+\vi A.$ Some assumptions are made about the nonlinearity $g$ such as $G^\p=g,$ for some $G\in C^1(H^1(\R^N;\C);\R).$ A solution is obtained as a non zero critical point of the real functional $F$ defined on $H^1_{A,V}(\R^N;\C)$ by,
\begin{gather*}
F(u)=\frac1{2\lambda}\|u\|_{H^1_{A,V}(\R^N)}^2 - G(u).
\end{gather*}
As a consequence, all the functional spaces are considered over the field of real numbers. It is shown that the equation makes sense in the topological dual space $H^{-1}_{A,V}(\R^N)$ of $H^1_{A,V}(\R^N;\C),$ and that the following holds.
\begin{gather*}
\Dr(\Omega;\C)\inj H^1_{A,V}(\R^N;\C)\inj H^1(\R^N;\C),		\\
H^{-1}(\R^N)\inj H^{-1}_{A,V}(\R^N)\inj \Dr^\p(\R^N),
\end{gather*}
with, in each case, dense embedding. Finally,
\begin{gather*}
-\Delta_Au=-\Delta u-\vi\nabla.(Au)-\vi A.\nabla u+|A|^2u, \; \text{ in } H^{-1}(\R^N),
\end{gather*}
and the product by $\vi$ is understood as in Definition~\ref{defiT}.

\subsection*{A nonlinear Schrödinger equation with external magnetic field}

In Schindler and Tintarev~\cite{MR1912753}, the authors consider the equation
\begin{gather*}
\begin{cases}
\left(\dfrac1\vi\nabla+A(x)\right)^2u+V(x)u=|u|^{q-1}u, \text{ in } \Omega\subseteq\R^3,	\medskip \\
u\in H^1_0(\Omega;\C)\setminus\{0\}.
\end{cases}
\end{gather*}
They show that a solution exists if $\Omega$ is asymptotically contractive (Schindler and Tintarev~\cite[Definition~2.6]{MR1912753}), $q\in(1,5),$ and $A\in L^2_\loc(\Omega;\R^3)$ and $V\in V^1_\loc(\Omega;\R)$ satisfy,
\begin{gather*}
\inf_{x\in\Omega}(V(x)+\lambda_0)>0,					\\
|A-\nabla\phi|<\liminf_{|y|\to\infty}V(y)-V, \text{ a.e.\:in } \Omega,
\end{gather*}
for some $\phi\in H^1_\loc(\Omega;\R),$  where $\lambda_0=\dsp\inf_{\left\{\substack{u\in\Dr(\Omega;\R) \hfill \\ \|u\|_{L^2(\Omega)}=1}\right.}\vint_\Omega|\nabla u(x)|^2\d x.$ To this end, they prove the existence of a minimizer of the real functional,
\begin{gather*}
\inf_{\left\{\substack{u\in H^1_0(\Omega;\C) \hfill \\ \|u\|_{L^{q+1}(\Omega)}=1}\right.}I(u),
\end{gather*}
where
\begin{gather*}
I(u)=\vint_\Omega\left(\left|\left(\dfrac1\vi\nabla+A(x)\right)u(x)\right|^2+V(x)|u(x)|^2\right)\d x.
\end{gather*}
Then the differentiability of $I$ is needed and all the functional spaces have to be considered over the field of real numbers. Here again, the product by $\vi$ is understood as in Definition~\ref{defiT}.

\subsection*{Existence in the nonlinear Schrödinger equation with bounded magnetic field}

Let $\nabla_A=\nabla+\vi A$ and $H^{1,2}_A(\R^N;\C)=\dot H^{1,2}_A(\R^N;\C)\cap L^2(\R^N;\C),$ where $\dot H^{1,2}_A(\R^N;\C)$ is the completion of $\Dr(\R^N;\C)$ with respect to the norm $\|\nabla_Au\|_{L^2(\R^N)}.$ In Schindler and Tintarev~\cite{MR4404858}, the authors prove that the minimization problem
\begin{gather*}
\inf_{\left\{\substack{u\in H^{1,2}_A(\R^N;\C) \hfill \\ \|u\|_{L^p(\Omega)}=1}\right.}J_{A,V}(u),
\end{gather*}
where $N\ge2,$ $p\in\left(2,\frac{2N}{N-2}\right)$ and 
\begin{gather*}
J_{A,V}(u)=\vint_{\R^N}\left(|\nabla_A u(x)|^2+V(x)|u(x)|^2\right)\d x,
\end{gather*}
admits a solution $u_\star,$ under some suitable assumptions satisfied by $A\in C^1(\R^N;\R^N)$ and $V\in C(\R^N;\R).$ Dealing with real Banach spaces to have that $J_{A,V}$ is differentiable, they obtain that $u_\star$ solve,
\begin{gather*}
\begin{cases}
-(\nabla_Au)^2+V(x)u=|u|^{p-1}u \; \text{ in } \R^N,	\medskip \\
u\in H^1_{A,V}(\R^N)\setminus\{0\}.
\end{cases}
\end{gather*}
Here again, the product by $\vi$ is understood as in Definition~\ref{defiT}.

\section*{Acknowledgements}
\baselineskip .5cm
The author acknowledges funding from ANR under grant ANR-17-EURE-0010 (Investissements d’Ave\-nir program). He is also grateful to Thierry Cazenave for interesting discussions during the preparation of this paper.

\baselineskip .4cm


\end{document}